\documentclass[journal]{IEEEtran}
\ifCLASSINFOpdf
\else
\fi
\usepackage{bm,algorithmicx,algpseudocode,threeparttable,booktabs}
\usepackage{amsmath,color}
\usepackage{cases}
\usepackage{stmaryrd}
\usepackage{amsfonts}
\usepackage{graphicx,times,psfig,amsmath}
\usepackage{cite}% Add all your packages here
\usepackage{hyperref}
\usepackage{multirow}
\usepackage{subfigure}

\begin{document}
\title{Multiagent based state transition algorithm for global optimization}

\author{Xiaojun~Zhou% <-this % stops a space
\thanks{X. Zhou is with the School
	of Automation, Central South University, Changsha, Hunan, China e-mail: (michael.x.zhou@csu.edu.cn) }% <-this % stops a space
\thanks{Manuscript received **, 2021; revised **.}}


\markboth{Journal of \LaTeX\ Class Files,~Vol.~6, No.~1, January~2007}%
{Shell \MakeLowercase{\textit{et al.}}: Bare Demo of IEEEtran.cls for Journals}



\maketitle

\begin{abstract}
%\boldmath
In this paper, a novel multiagent based state transition optimization algorithm with linear convergence rate named MASTA is constructed. It first generates an initial population randomly and uniformly. Then, it applies the basic state transition algorithm (STA) to the population and generates a new population. After that, It computes the fitness values of all individuals and finds the best individuals in the new population. Moreover, it performs an effective communication operation and updates the population. With the above
iterative process, the best optimal solution is found out. Experimental results based on some common benchmark functions and comparison with some stat-of-the-art optimization algorithms, the proposed MASTA algorithm has shown very superior and comparable performance.
\end{abstract}


% \fi
\IEEEpeerreviewmaketitle

\section{Introduction}
\IEEEPARstart{S}ince Professor John Holland from the University of Michigan proposed the genetic algorithm (GA) in the early 1970s \cite{Holland1984}, the intelligent optimization algorithms represented by GA \cite{Goldberg1989} have made considerable development, and many other new intelligent optimization algorithms have been sprung up, such as the simulated annealing, ant colony algorithm and particle swarm optimization. Research on the intelligent optimization algorithms has been becoming the most active research direction in the intelligence science, information science and artificial intelligence, and got rapid popularization and application in many engineering fields \cite{Yang2010}. Most present intelligent optimization algorithms are based on behaviorism and imitation learning, and they solve complex optimization problems by simulating the biological evolution of birds, bees, fish and other organisms in nature \cite{karaboga2007powerful,karaboga2014quick,kulkarni2011particle}. However, the behaviorism based intelligent optimization algorithms mainly focus on the imitation, and they imitate and learn what they encounter \cite{kirkpatrick1983optimization}. It is too mechanical and has great blindness, and cannot profoundly reflect the essence, purpose and requirements of the optimization algorithms \cite{karafotias2015parameter}. On the one hand, these appearance imitation learning based algorithms have resulted in the poor scalability. Most intelligent optimization algorithms perform well on some problems with low dimensions, but show significantly deteriorated performance on problems the with higher dimensions. On the other hand, it makes the algorithms to prone to weird phenomena such as the stagnation or premature convergence, that is, the algorithms may be stuck at any random points, rather than the optimal solution mathematically \cite{mallipeddi2011differential,mitra1985convergence}. In order to eliminate the stagnation, improve the scalability and expand the application scope of the intelligent optimization algorithms, Dr. Zhou has proposed a new intelligent optimization algorithm based on the structured learning called state transfer algorithm (STA) in 2012 \cite{zhou2011initial,zhou2011new,zhou2012state,zhou2018statistical,zhou2019dynamic,zhou2018hybrid}.

STA is an intelligent stochastic global optimization algorithm based on structured learning. It well reflects the essence, purpose and requirements of the optimization algorithms and has the system framework consisting of five core structural elements including globality, optimality, rapidity, convergence and controllability.
The initial version of the state transition algorithm (STA) algorithm \cite{zhou2011initial} was proposed on the basis of the concepts of state and state transition. STA introduced the three special operators named rotation, translation and expansion for continuous function optimization problems and general elementary transformation operator for discrete function optimization problems, and illustrated very good search capability based on a discrete problem and four common benchmark continuous functions.
In order to promote the global search ability of STA, \cite{zhou2011new} designed an axesion operator to search along the axes and strengthen single dimensional search.
\cite{zhou2012state} studied the adjusting measures of the transformations to keep the balance of exploration and exploitation and discussed the convergence analysis about STA based on random search theory, and introduced the communication strategy into the basic algorithm and presented intermittent exchange to prevent premature convergence.
\cite{zhou2018hybrid} presented an efficient hybrid feature selection method based on binary state transition algorithm and ReliefF containing the filter and wrapper two phases to solve the feature selection problems.
\cite{zhou2014constrained} presented a population-based continuous state transition algorithm named constrained STA to solve the continuous constrained optimization problems, and the constrained STA used a two-stage strategy that in the early stage of an iteration process the feasibility preference method is adopted and in the later stage it is changed to the penalty function method. Based on several benchmark tests the performance of the constrained STA with a two-stage strategy were proved to outperform other single strategy in terms of both global search ability and solution precision.
\cite{zhou2014nonlinear} introduce the STA to transform identification and control for nonlinear system into optimization problems. It first applied STA to identify the optimal parameters of the estimated system with previously known structure and then designed an off-line PID controller optimally by using the STA.
\cite{zhou2016comparative} conducted a comparative study of the performance of standard continuous STA with the comparison with some other state-of-the-art evolutionary algorithms on large scale global optimization problems and experimental results showed that its global search ability is much superior to the competitors.
In order to solve unconstrained integer optimization problems, \cite{zhou2016discrete} presented a discrete state transition algorithm, using several intelligent operators for local exploitation and global exploration and a dynamic adjustment strategy to capture global solutions with high probability.
\cite{zhou2016optimal} introduced a discrete STA to solve several cases of water distribution networks with a parametric study of the `restoration probability and risk probability' in the dynamic STA and the investigation of the influence of the penalty coefficient and search enforcement on the performance of the algorithm.
\cite{zhou2018setpoint} proposed a multi-objective state transition algorithm (MOSTA) for optimizing the parameters of the PID controllers.
\cite{zhou2018dynamic} proposed a stochastic intelligent optimization method based on STA in which a novel dynamic adjustment strategy called `risk and restoration in probability' is incorporated into STA to transcend local optimality, and designed a refined dynamic state transition algorithm with the appropriately chosen risk probability and restoration probability, in order to solve the sensor network localization problem without additional assumptions and conditions on the problem structure.

Generally speaking, the basic idea of STA is to consider each solution of the optimization problem as a state, and the iterative updating process of the solution as a state transition process. The discrete time state space expression in modern control theory is used as a unified framework for generating candidate solutions, and state transformation operators are designed based on the framework.
Unlike most population based evolutionary algorithms, the standard state transition algorithm is a kind of individual based evolutionary algorithm. It is based on the given current solutions, performs multiple independent runs of a certain state transformation operator to generate the candidate solution set by sampling ways, compares with the current solutions, and iteratively updates the current solutions until some termination conditions are met. It is worth mentioning that in the state transition algorithm each state transformation operator can generate geometric neighborhood with the regular shape and controlled size. It designs different state transformation operators including the rotation, translation, expansion, and axesion transformation operators to satisfy some function requirements in the global search, local search and heuristic search, and uses different operators timely in the form of alternate rotation, in order to make the state transition algorithm to quickly find the global optimal solutions in probability. In recent few years, STA has been emerging as a novel intelligent optimization method for global optimization. Unlike the majority of population-based evolutionary algorithms, the basic STA is an individual-based optimization method. With the aim to find a possibly global optimum as soon as possible, in basic STA, it contains three main types of state transformation operators: global, local and heuristic operators. Global operator aims to generate a neighborhood, formed by all possible candidate solutions, which may contain the global optimum. Local operator aims to generate a neighborhood strictly in geometry, i.e., a hypersphere with the incumbent best solution as the center. And heuristic operators aims to find potentially useful solutions, to accelerate the convergence speed and to avoid blind enumeration. To be more specific, with four transformation operators in basic continuous STA, the rotation transformation is designed for local search, the expansion transformation is designed for global search, and the translation and axesion transformation are proposed for heuristic search. More importantly, the size of the neighborhood generated by any state transformation operator can be controlled with a parameter. Together with the sampling technique and alternative utilization, the basic STA has exhibited powerful search ability in global optimization.

To further accelerate the convergence of state transition algorithm, in this study, we aims to propose a multiagent based state transition optimization algorithm with adequate convergence rate.

\section{Preliminaries}
Let us focus on the following unconstrained continuous optimization problem
\begin{equation}
\min_{x \in \mathbb{R}^n} f(x)
\end{equation}
where $f(x)$ is a continuous function, and $n$ is the size the decision vector $x$.

\subsection{Basic state transition algorithm}
The basic state transition algorithm (STA) is an individual-based optimization method \cite{zhou2011initial,zhou2011new}.
On the basis of state space representation, the unified form of
solution generation in basic STA can be formulated as follows
\begin{equation}
x_{k+1}= A_{k} x_{k} + B_{k} u_{k}
\end{equation}
where $x_{k}$ and $x_{k+1}$  stand for a current state and the next state respectively, corresponding to candidate solutions of the optimization problem; $u_{k}$ is a function of $x_{k}$ and historical states;
$A_{k}$ and $B_{k}$ are state transition matrices (sometimes with random entries), integrating with $x_{k}$ and $u_{k}$, forming the transformation operators.

In basic continuous STA, it contains four different state transformation operators, including rotation, translation, expansion and axesion transformation respectively, as formulated as follows
\begin{eqnarray}
x_{k+1} &=& x_{k}+\alpha \frac{1}{n \|x_{k}\|_{2}} R_{r} x_{k}, \nonumber \\
x_{k+1} &=& x_{k}+  \beta  R_{t}  \frac{x_{k}- x_{k-1}}{\|x_{k}- x_{k-1}\|_{2}}, \nonumber \\
x_{k+1} &=& x_{k}+  \gamma  R_{e} x_{k},\nonumber \\
x_{k+1} &=& x_{k}+  \delta  R_{a}  x_{k}
\end{eqnarray}
where $\alpha, \beta, \gamma, \delta$ are state transformation factors, and $R_{r}$ $\in$ $\mathbb{R}^{n\times n}, R_{t}$ $\in \mathbb{R}, R_{e} \in \mathbb{R}^{n \times n}, R_{a}$ $\in \mathbb{R}^{n \times n}$ are different random matrices.
The details of the basic continuous STA can be referred to \cite{zhou2011initial}.

\newtheorem{example}{Example}
\newtheorem{proposition}{Proposition}
\newtheorem{remark}{Remark}

\subsection{Rate of linear convergence}
In numerical optimization, the linear convergence of an optimization algorithm can be defined as follows
\begin{equation}
\frac{\|x_{k+1} - x^{*}\|}{\|x_{k} - x^{*}\|} = \eta \;\;\;\;
\end{equation}
where $x_{k} \in \mathbb{R}^n$ is the $k$th iterative point, $x^{*} \in \mathbb{R}^n$ is a nominal optimal point, and $\eta$ (usually $0 < \eta < 1$) is called the rate of convergence.

\section{Multiagent based state transition optimization algorithms}
In this section, we aim to construct effective multiagent-based state transition optimization algorithms.
Supposing that there exist a population with $N$ agents (individuals) and each agent has the same (similar) dynamics, inspired by the concept of linear convergence and multiagent system,
a linear form of solution generation in multiagent-based STA for the $i$the agent is given by
\begin{equation}
x^{i}_{k+1}= A_{k} x^{i}_{k} + B_{k} u^{i}_{k}
\end{equation}
where $u^{i}_{k}$ is a function of current state and historical states of all agents.

\subsection{Multiagent based STAs with guaranteed convergence}

\begin{figure}[!htbp]
\centering
  \includegraphics[width=8cm]{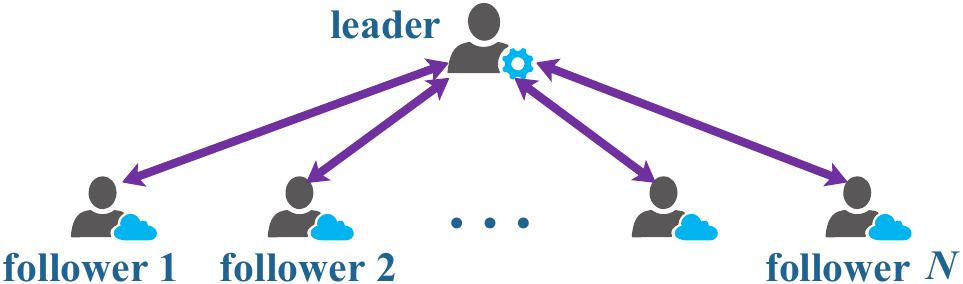}
\caption{Illustration of the leader-follower multiagent systems of type I}
\label{fig:leader-follower_typeI}       % Give a unique label
\end{figure}

Considering the leader-follower multiagent systems of type I as illustrated in Fig. \ref{fig:leader-follower_typeI}, where the
leader represents the best individual (with the smallest function value), which is selected from the population, and the $i$th follower can only communicate information with the leader. From a numerical optimization viewpoint, the leader is selected from the followers can be changed during the iterative process.

\subsubsection{Multiagent-based STA with a fixed rate of convergence}
Firstly, the following iterative equation for multiagent-based state transition optimizer for the $i$th agent is constructed
\begin{equation}
\label{eq_fixed_rate}
x^{i}_{k+1} = \eta x^{i}_{k} + (1-\eta)Best
\end{equation}
where $\eta (|\eta| < 1)$ is a constant, $x^i_{k}$ is the $i$th follower and $Best$ is the leader (incumbent best solution) in the population.

From the perspective of convergence theory, if the incumbent best solution $Best$ remains unchanged, then we have
\begin{equation}
x^{i}_{k} = Best + \eta^{k} (x^{i}_{0} - Best)
\end{equation}
that is to say, if $|\eta| < 1$, each agent $x^{i}$ converges  to $Best$,
\begin{equation}
\lim_{k \rightarrow \infty} x^{i}_{k} =  Best  + \lim_{k \rightarrow \infty}  [\eta^{k} (x^{i}_{0} - Best)] = Best
\end{equation}

\begin{example}
\label{exam_fixed_rate}
Considering the following unconstrained optimization problem
\begin{equation*}
f_0(x) = (x_1 - 1)^2 + (x_2 - 2x_1^2)^2
\end{equation*}
By using the Eq.(\ref{eq_fixed_rate}) with a fixed rate ($\eta = -0.5$) well distributed initial population,
trajectories of all followers and the leader are given in Fig. \ref{fig:trajectories_fixed_rate}.
It is found that all agents can converge to the optimal solution.
\begin{figure*}
\centering
  \includegraphics[width=8cm]{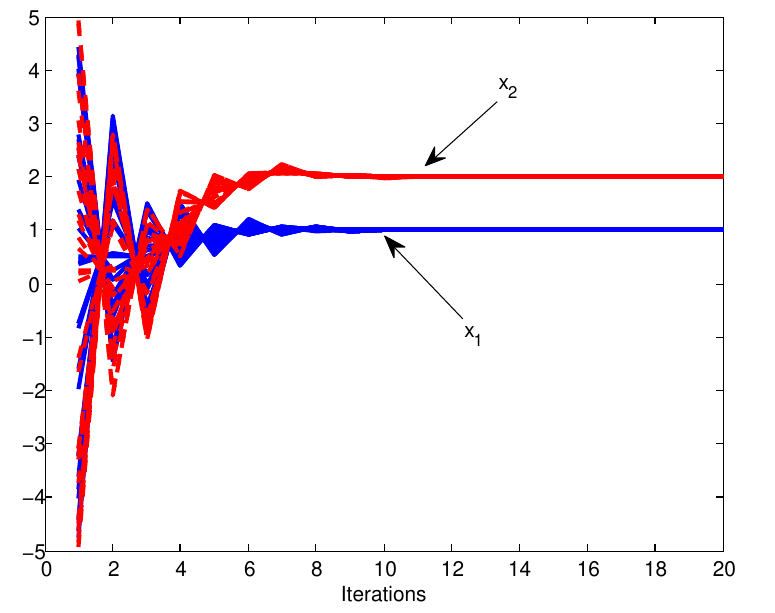}
  \includegraphics[width=8cm]{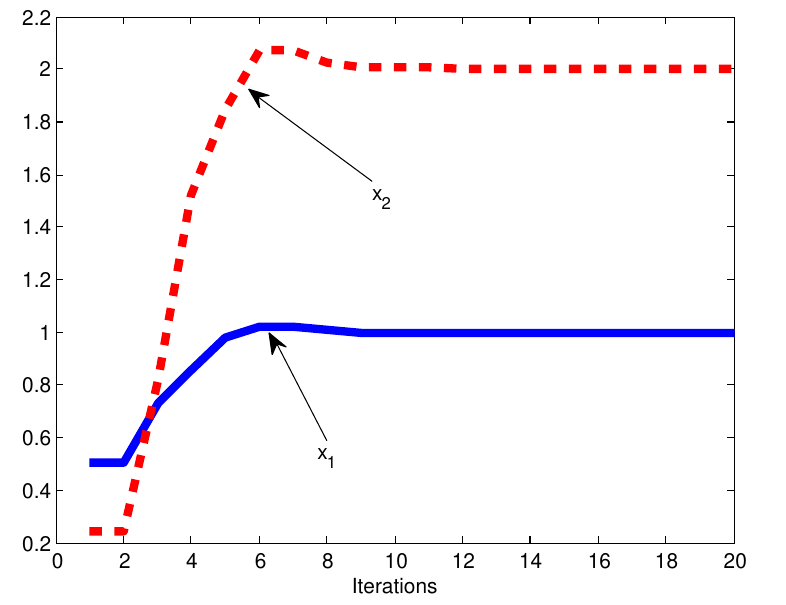}
\caption{Trajectories of all followers and the leader with a fixed rate.}
\label{fig:trajectories_fixed_rate}       % Give a unique label
\end{figure*}
\end{example}

\subsubsection{Multiagent-based STA with a varying rate of convergence}
Since it is hard to choose an appropriate rate of convergence $\eta$, let us consider the following iterative formula with a varying rate of convergence
\begin{equation}
\label{eq_varying_rate}
x^{i}_{k+1} = \eta_k x^{i}_{k} + (1 - \eta_k) Best
\end{equation}
Then, we can obtain
\begin{equation}
x^{i}_{k+1} = Best + \eta_0 \eta_1 \cdots \eta_k  (x^{i}_{0} - Best)
\end{equation}
If  $|\eta_k| < 1$, since the sequence $\{|\eta_0 \eta_1 \cdots \eta_k|\}_{k=0}^{\infty}$ is strictly decreasing, then each agent $x^{i}$ converges to $Best$
\begin{equation}
\lim_{k \rightarrow \infty} x^{i}_{k} =  Best + \lim_{k \rightarrow \infty} \eta_0 \eta_1 \cdots \eta_k  (x^{i}_{0} - Best) = Best
\end{equation}

\begin{example}
\label{exam_varying_rate}
Considering the same optimization problem as in Example \ref{exam_fixed_rate}, by using the Eq.(\ref{eq_varying_rate}) with a linearly increasing rate ($\eta: -0.9 \rightarrow -0.1$) and well distributed initial population, trajectories of all followers and the leader are given in Fig. \ref{fig:trajectories_varying_rate}. It is found that all agents can converge to the optimal solution as well.

\begin{figure*}[!htbp]
\centering
  \includegraphics[width=8cm]{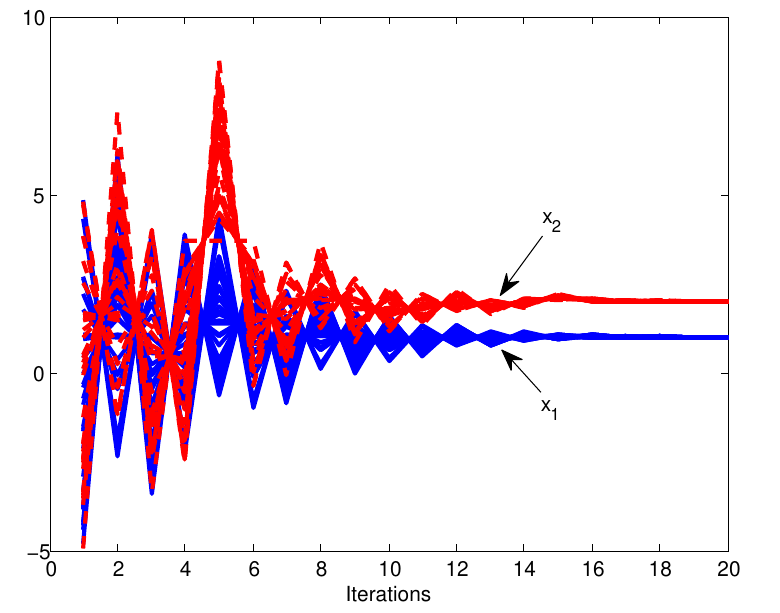}
  \includegraphics[width=8cm]{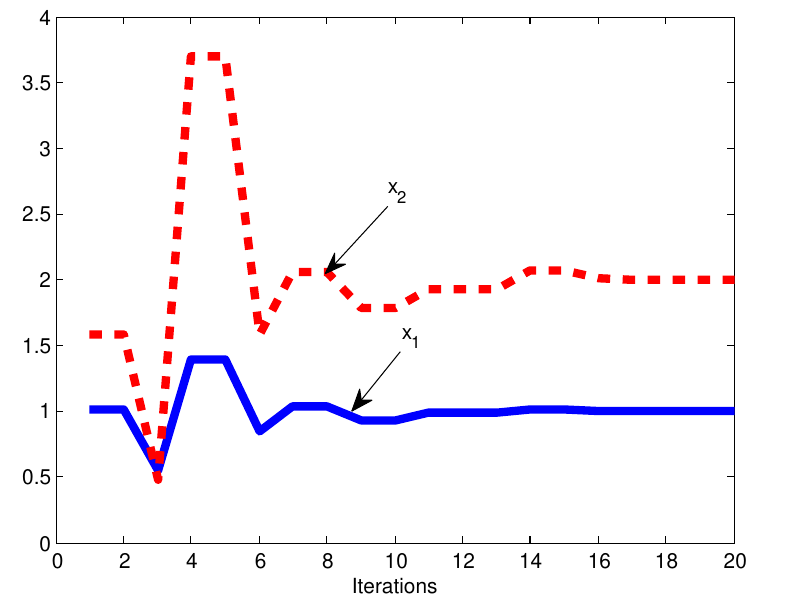}
\caption{Trajectories of all followers and the leader with a varying rate.}
\label{fig:trajectories_varying_rate}       % Give a unique label
\end{figure*}
\end{example}

\subsubsection{Multiagent-based STA with a stochastic rate of convergence}
Let us consider the following iterative formula with flexible rate of convergence
\begin{equation}
\label{flexible_rate}
x^{i}_{k+1} = \hat{\eta} x^{i}_{k} + (1 - \hat{\eta})Best
\end{equation}
where $\hat{\eta}$ is a random number.

If $\hat{\eta}$ is in the interval $(-1,1)$, it is not difficult to understand that the sequence $\{x^{i}\}_{k=0}^{\infty}$ generated by Eq.(\ref{flexible_rate})  converges to $Best$.

\begin{proposition}
If $\hat{\eta}$ is uniformly distributed random number in the interval $(-2,2)$, then the sequence $\{x^{i}\}_{k=0}^{\infty}$ generated by Eq.(\ref{flexible_rate}) converges to $Best$.
\end{proposition}

Furthermore, considering the following iterative formula with $L$ adjustable parameters
\begin{equation}
\label{flexible_rate_generalized}
x^{i}_{k+1} = Best + \hat{\eta}^1 \hat{\eta}^2 \cdots \hat{\eta}^L (x^{i}_{k} - Best)
\end{equation}

\begin{proposition}
If $\hat{\eta}^1, \hat{\eta}^2, \cdots, \hat{\eta}^L$ are uniformly distributed random number in the interval $(-2,2)$, then the sequence $\{x^{i}\}_{k=0}^{\infty}$ generated by Eq.(\ref{flexible_rate_generalized})  converges to $Best$.
\end{proposition}

\begin{example}
Considering the following iterative formula
\begin{equation} \label{equ:IterativeFormula}
x_{k+1} = \hat{\eta} x_{k} + (1 - \hat{\eta})Best
\end{equation}
where $\hat{\eta}$ is uniformly distributed random number in the interval $(-2,2)$. Let $Best = [1,2]^T$,
for any given initial point $x_0$, the trajectories using the uniformly distribution will converge to the $Best$, as illustrated in Fig.\ref{fig:trajectories_flexible}.
\begin{figure}[!htbp]
\centering
  \includegraphics[width=8cm]{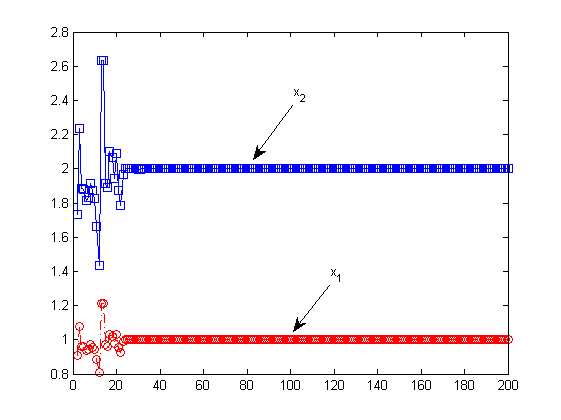}
\caption{Illustration of the trajectories by the uniformly distribution}
\label{fig:trajectories_flexible}       % Give a unique label
\end{figure}
\end{example}

\begin{proposition}
If $\hat{\eta}^1, \hat{\eta}^2, \cdots, \hat{\eta}^L$ are Gaussian distributed random numbers in the interval $(-2,2)$, then the sequence $\{x^{i}\}_{k=0}^{\infty}$ generated by Eq.(\ref{flexible_rate_generalized})  converges to $Best$.
\end{proposition}

\begin{figure}[!htbp]
\centering
  \includegraphics[width=8cm]{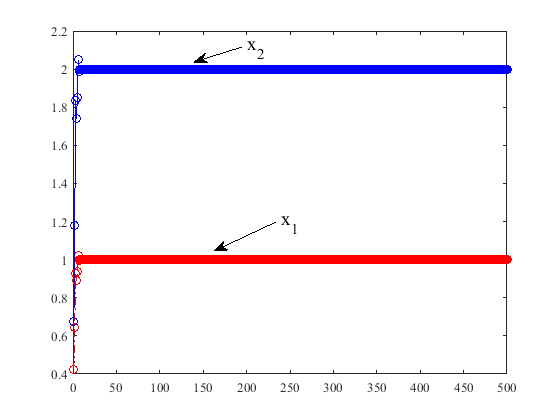}
\caption{Illustration of the trajectories by the Gaussian distribution}
\label{fig:gaussian_distribution_trajectories_flexible}       % Give a unique label
\end{figure}

Take the iterative formula by Eq. (\ref{equ:IterativeFormula}) where $\hat{\eta}$ is Gaussian distributed random number in the interval $(-2,2)$ for an example. Let $Best = [1,2]^T$, for any given initial point $x_0$, the trajectories using the Gaussian distribution will converge to the $Best$, as illustrated in Fig.\ref{fig:gaussian_distribution_trajectories_flexible}.

\subsection{Multiagent-based STAs with guaranteed optimality in probability}
\subsubsection{Multiagent-based STA with symmetry operation}
As discussed above, if a multiagent system \cite{Dorri2018} is constructed based on any of the mentioned three types: fixed rate of convergence, changeable rate of convergence, or flexible rate of convergence, the resulting multiagent-based STA is convergent. However, there is no guarantee that those proposed algorithms can converge to an optimum.

\begin{figure}[!htbp]
\centering
  \includegraphics{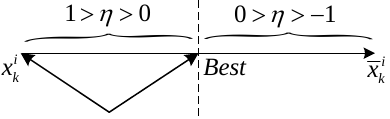}
\caption{Illustration of the fixed rate of convergence: $x^{i}_{k+1} = \eta x^{i}_{k} + (1-\eta)Best$}
\label{fig:fixed-rate}       % Give a unique label
\end{figure}

To explain this point, let us take the fixed rate of convergence in Eq. (\ref{eq_fixed_rate}) for exmaple. For the $i$th individual,
as shown in Fig. \ref{fig:fixed-rate}, if $1 > \eta > 0$, the next iterative $x^i_{k+1}$ will be on
the line segment that goes through from $x^i_{k}$ to $Best$. If $0 > \eta > -1$, the next iterative $x^i_{k+1}$ will be on the line segment that goes through from $Best$ to $\overline{x}^i_{k}$ (the symmetric point of $x^i_{k}$, i.e., $\overline{x}^i_{k} - Best = Best - x^i_{k})$.

Considering the $k$th iteration with multiple agents ($x^{i}_{k}, i = 1, 2, \cdots$, $N$), to make sure that a multiagent-based STA can converge to an optimum in probability at the next iteration, i.e., $Best$ is an optimum in probability, the following four conditions should be satisfied simultaneously:
\begin{itemize}
  \item[(1)] $|\eta|$ is sufficient small;
  \item[(2)] $N$ is sufficient large;
  \item[(3)] all candidate solutions $x^{i}_{k+1}$ should be distributed uniformly around $Best$;
  \item[(4)] $f(Best) < f(x^{i}_{k+1})$, $\forall i = 1, 2, \cdots, N$.
\end{itemize}

\begin{figure*}[!htbp]
\centering
  \includegraphics[width=5.5cm]{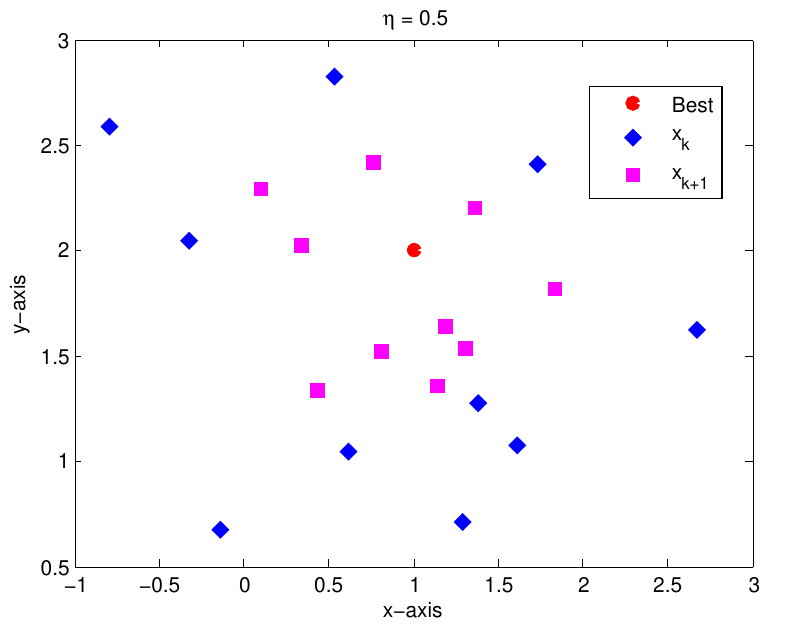}
  \includegraphics[width=5.5cm]{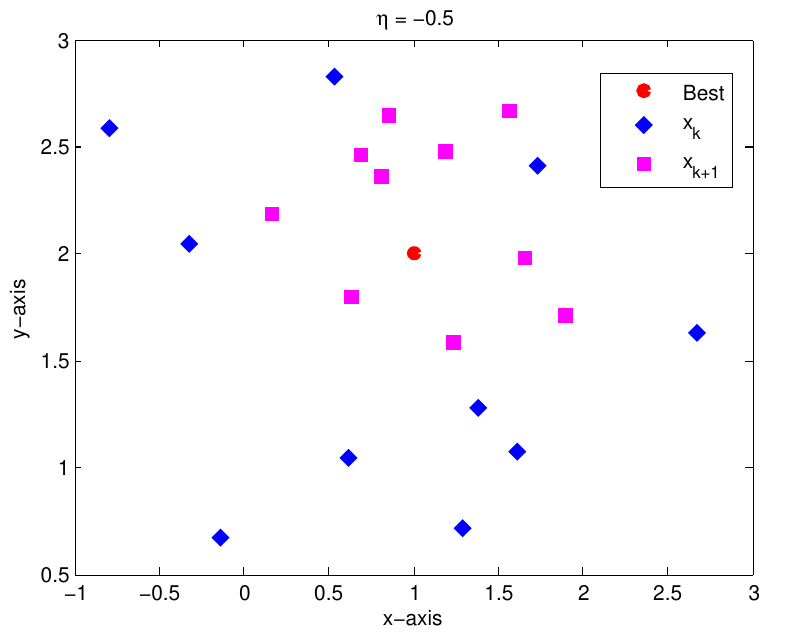}
\caption{Illustration of next population if current population is distributed uniformly around $Best$.}
\label{fig:distributed-uniformly}       % Give a unique label
\end{figure*}

\begin{figure*}%[!htbp]
\centering
  \includegraphics[width=5.5cm]{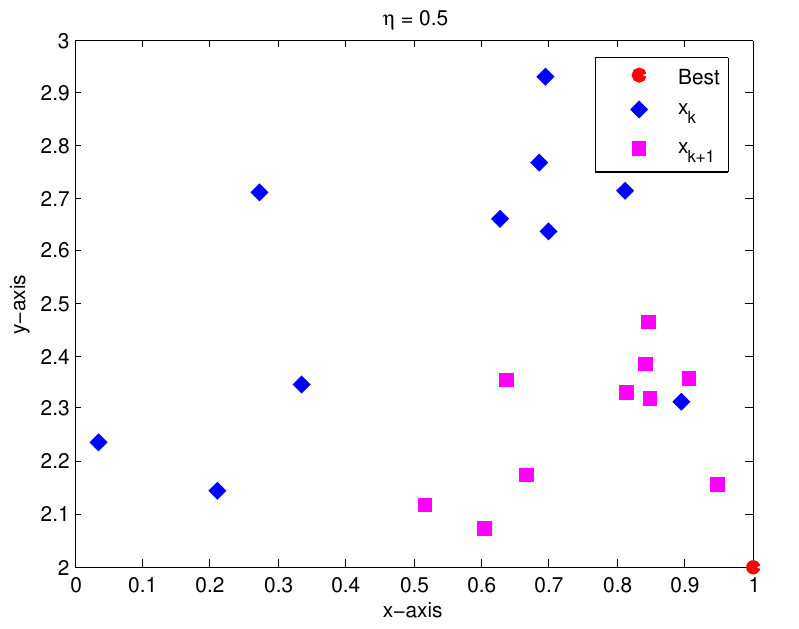}
  \includegraphics[width=5.5cm]{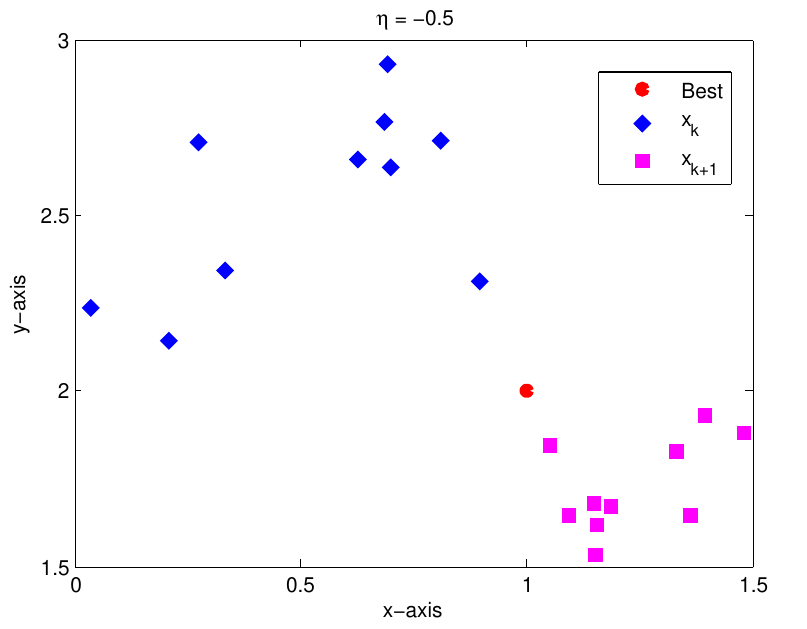}
\caption{Illustration of next population if current population is not distributed uniformly around $Best$.}
\label{fig:distributed-nonuniformly}       % Give a unique label
\end{figure*}

Among the four conditions, the third one is not easy to meet. As illustrated from Fig. \ref{fig:distributed-uniformly}
to Fig. \ref{fig:distributed-nonuniformly}, if the current population is distributed uniformly around $Best$, the next population is most likely to be distributed uniformly around $Best$ as well, but that is not true vice versa.

Since the incumbent best solution $Best$ may be changed if a better solution is found in the next population, it is very likely that the next population is not distributed uniformly around the nest best solution. As a result, the next population should be repaired to meet the third condition.
Considering that $2Best - x^{i}_{k+1} = Best - \eta (x^{i}_{k} - Best)$, the following symmetry operation is used:
\begin{equation}
\overline{x}^{i}_{k+1} = 2Best - x^{i}_{k+1}.
\end{equation}
where $x^{i}_{k+1}$ can be generated by any of the three types of convergence.

Combing both the candidate solutions $x^{i}_{k+1}$ and $\overline{x}^{i}_{k+1}$ by using the symmetry operation, together with the flexible rate of convergence,
it is not difficult to imagine that the repaired population is much more likely to be distributed uniformly around the incumbent best solution.

\subsubsection{Multiagent-based STA with convex combination}
To generate a better distribution for the next iteration, the following new iterative equation for the $i$th
agent is constructed
\begin{equation}
\label{eq_convex_combination}
x^i_{k+1} = \eta x^i_{k} +  \zeta x^j_{k} + (1-\eta-\zeta) Best
\end{equation}
where $0 \leq \eta, 0 \leq \zeta$ and $\eta + \zeta \leq 1$, and they can be constant, varying or random.
The $j$th agent is randomly selected from the rest of the population.

\begin{figure}[!htbp]
\centering
  \includegraphics[width=8cm]{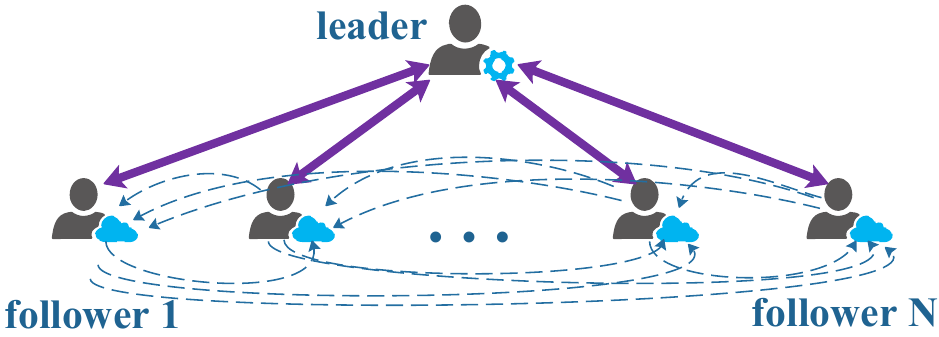}
\caption{Illustration of the leader-follower multiagent systems of type II}
\label{fig:leader-follower_typeII}       % Give a unique label
\end{figure}

\begin{figure}[!htbp]
\centering
  \includegraphics{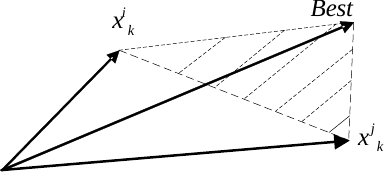}
\caption{Illustration of the convex combination}
\label{fig:convex_combination}       % Give a unique label
\end{figure}

Actually, Eq. \ref{eq_convex_combination} is a convex combination of three points ($x^i_{k}, x^j_{k}$ and $Best$), and the resulting candidate solution $x^i_{k+1}$ will lie on a triangular plane, as illustrated in Fig. \ref{fig:convex_combination}.

\begin{figure*}[!htbp]
\centering
  \includegraphics[clip=true, width=8cm]{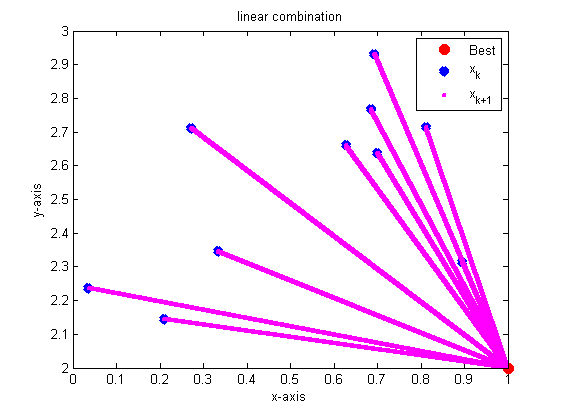}
  \includegraphics[clip=true, width=8cm]{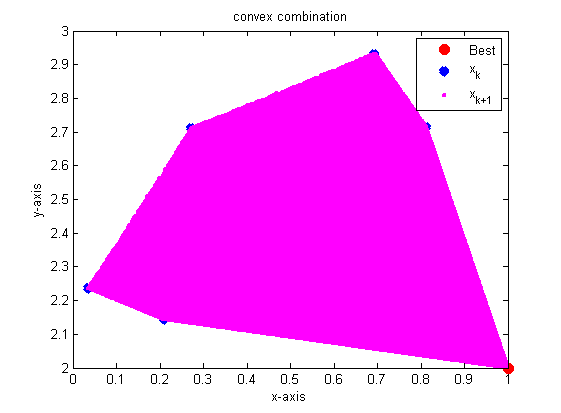}
\caption{Illustration of the linear and convex combination with flexible rate}
\label{fig:linear_convex_combination}       % Give a unique label
\end{figure*}

With the same current population as in Fig. \ref{fig:distributed-nonuniformly}, the distribution of the next population using linear combination Eq. (\ref{flexible_rate}) and convex combination Eq.(\ref{eq_convex_combination}) with flexible rate is illustrated in Fig. \ref{fig:linear_convex_combination}. It is found the distribution of convex combination is much wider than that of linear combination.

\subsubsection{Multiagent-based STA with elementwise operation}
\begin{figure}[!htbp]
\centering
  \includegraphics[clip=true, width=8cm]{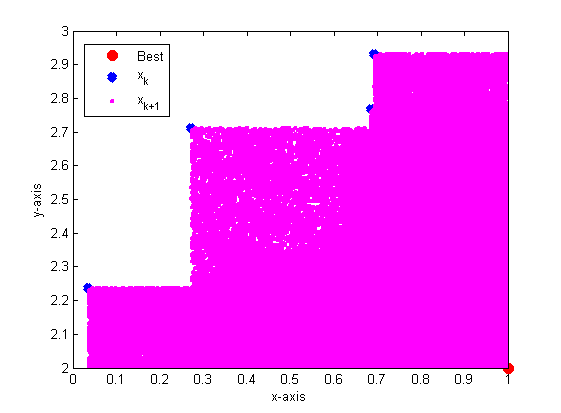}
\caption{Illustration of the linear combination with flexible rate and elementwise operation}
\label{fig:linear_elementwise}       % Give a unique label
\end{figure}

Although the distribution generated by the classic linear combination is not so good, the elementwise operation can be used to make it much more useful, as given below
\begin{equation}
\label{eq_linear_elementwise}
x^{i}_{k+1,d} = Best_{k,d} + \eta_{k,d} (x^{i}_{k,d} - Best_{k,d})
\end{equation}
where $x^{i}_{k,d}$ means the $d$th component of $i$th agent for the $k$th iteration, and $\eta_{k,d}$ is the rate of convergence, which can be constant, varying or random. That is to say, different rates of convergence are given for different dimensions.
With the same current population as in Fig. \ref{fig:distributed-nonuniformly}, the distribution of the next population using linear combination and elementwise operation Eq. (\ref{eq_linear_elementwise}) with flexible rate is illustrated in Fig. \ref{fig:linear_elementwise}.

\begin{remark}
The elementwise operation can also be extended to the convex combination, as shown below
\begin{equation}
\label{eq_elementwise_convex}
x^{i}_{k+1,d} = Best_{k,d} + \eta_{k,d} (x^{i}_{k,d} - Best_{k,d}) + \zeta_{k,d} (x^{j}_{k,d} - Best_{k,d})
\end{equation}
\end{remark}

\subsubsection{Multiagent-based STA with rotation transformation}
\begin{figure}[!htbp]
\centering
  \includegraphics[clip=true, width=8cm]{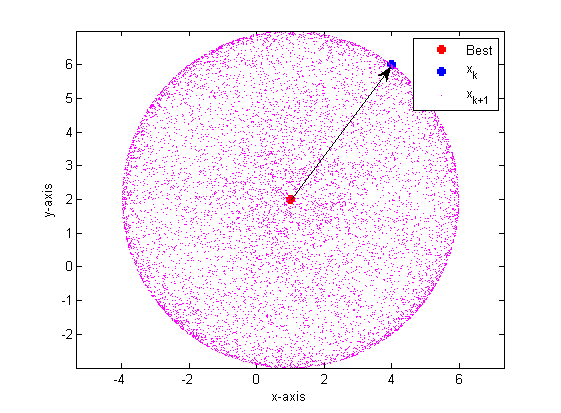}
\caption{Illustration of the new rotation transformation}
\label{fig:rotation}       % Give a unique label
\end{figure}

To further improve the distribution of the next population, a new rotation transformation is proposed as follows
\begin{equation}
x^i_{k+1} = Best + \eta \frac{R_r}{\|R_r\|} (x^i_{k} - Best)
\end{equation}
where $0< \eta \leq 1$ is a positive constant  and $R_r$ is a random matrix, usually with its entries distributed uniformly in the interval (-1,1). It is easy to prove that
\begin{eqnarray}
\|x^i_{k+1} - Best\| \leq  \eta  \|x^i_{k} - Best\|
\end{eqnarray}
that is to say, the possible candidate solutions will lie in a hypersphere with $Best$ as the center and $\eta  \|x^i_{k} - Best\|$ as the radius, as illustrated in Fig. \ref{fig:rotation}.

\subsection{Efficient multiagent-based STA for global optimization}
As discussed above, with different types of rate of convergence and various strategies to improve the distribution of next population, many multiagent based state transition optimization algorithms can be proposed.
In this part, for simplicity, a multiagent based STA integrating the flexible rate of convergence with the new rotation transformation is focused.

\subsubsection{Alleviating premature convergence}
When the incumbent best solution keeps unchanged after the predefined number of iterations, the value of the parameter $\alpha$ decreases gradually. Thus, the smaller area surrounding the incumbent best solution is searched, and it is very easy to find out a much better optimal solution than the incumbent best solution. Therefore, the premature convergence is well alleviated and it is very promising to get the global optimal solution from the search region.

\subsubsection{Stopping criterion}
The search process is terminated when no improvement of the found best solution is observed from the search space during the search.

\subsection{The proposed multiagent-based STA algorithm}
\begin{figure}[!htbp]
\centering
  \includegraphics{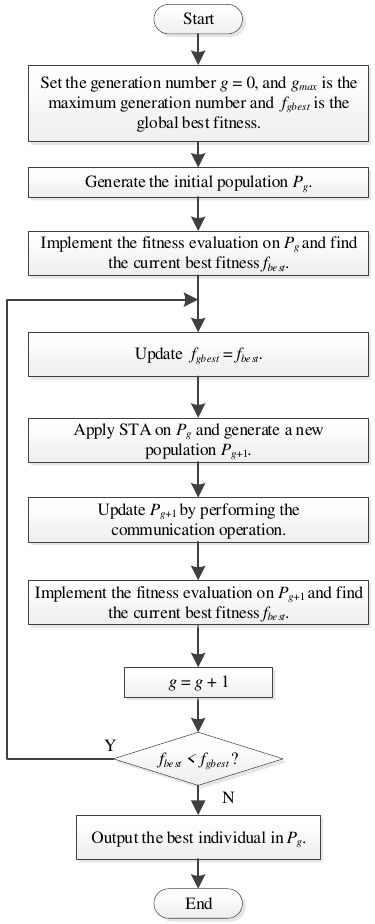}
\caption{The flowchart of the MASTA algorithm.}
\label{fig:flowchartalgorithm}       % Give a unique label
\end{figure}

In the paper, we have proposed a multiagent-based STA algorithm named MASTA. The MASTA algorithm carries on an iterative procedure as follows. It first generates an initial population $P_g$ containing $N$ solutions randomly and uniformly with the generation number $g=0$. Then, it updates the global best solution $x_{gbest}$ from the population $P_g$. After that, it applies the STA to the population $P_g$ and generates a new population $P_{g+1}$. Moreover, it performs the communication operation and updates the population $P_{g+1}$ for $CF$ iterations in which $CF$ is the communication frequency. In addition, it computes the fitness values of all solutions in the population $P_{g+1}$, and updates the global best solution $x_{gbest}$ and the generation number $g=g+1$.
If the stopping criterion, that is, the global best solution is not improved, has not been satisfied, the above process is repeated; otherwise, the best optimal solution is output from the population. The flowchart of the proposed MASTA algorithm is shown in Fig. \ref{fig:flowchartalgorithm}.

\begin{table*}[htbp]
\centering
\caption{Benchmark functions} \label{tab:BenchmarkFunctions}
    \begin{tabular}{cccc}
    \toprule
    Name of function& Function definition&  Range& $f_{min}$\\
    \midrule
    Spherical       &${f_1} = \sum\nolimits_{i = 1}^n {x_i^2} $	    &[-100,100]	    &0\\
    Rastrigin       &${f_{\rm{2}}} = \sum\nolimits_{i = 1}^n {\left( {x_i^2{\rm{ - 10}}\cos (2\pi {x_i}) + 10} \right)} $ &[-5.12,5.12] & 0\\
    Griewank        &${f_3} = \frac{1}{{4000}}\sum\nolimits_{i = 1}^n {x_i^2}  - \prod\nolimits_i^n {\cos \left| {\frac{{{x_i}}}{{\sqrt i }}} \right|}  + 1$ &[-600,600] & 0\\
    Rosenbrock      &${f_4} = \sum\nolimits_{i = 1}^n {\left( {100{{\left( {{x_{i + 1}}{\rm{ - }}x_i^2} \right)}^2} + {{\left( {{x_i} - 1} \right)}^2}} \right)} $ &[-30,30] &0\\
    Ackley          &${f_{\rm{5}}} = {\rm{20 + e - 20exp}}\left( { - 0.2\sqrt {\frac{1}{n}\sum\nolimits_{i = 1}^n {x_i^2} } } \right)$	&[-32,32]	&0\\
    Schaffer        &${f_6} = {\rm{0}}{\rm{.5 + }}\frac{{\sin {{\left( {\sqrt {x_1^2 + x_2^2} } \right)}^2} - 0.5}}{{{{\left( {1 + 0.001(x_1^2 + x_2^2)} \right)}^2}}}$ &[-100,100] &0 \\
    Easom           &${f_7} =  - \cos ({x_1})\cos ({x_2})\exp ( - ({x_1} - \pi ) - ({x_2} - \pi ))$ &[-100,100] &-1 \\
    Goldstein-Price &$\begin{array}{l}
    {f_8} = [1 + {\left( {{x_1} + {x_2} + 1} \right)^2}(19 - 14{x_1} + 3x_1^2\\
     - 14{x_2} + 6{x_1}{x_2} + 3x_2^2)] + [30 + {\left( {2{x_1} - 3{x_2}} \right)^2}\\
    (18 - 32{x_1} + 12x_1^2 + 48{x_2} - 36{x_1}{x_2} + 27x_2^2)]
    \end{array}$ &[-2,2] &3 \\
    \bottomrule
    \end{tabular}
\end{table*}

\section{Experimental study}
To test the performance of the proposed MASTA algorithm, we carried on several experiments on two dimensional, ten dimensional, thirty dimensional and fifty dimensional functions respectively.

\begin{table*}[htbp]
\centering
\caption{Experimental results on test functions(2D)} \label{tab:results2D}
    \begin{tabular}{cccccc}
    \toprule
    Fcn                           &Statistic &MASTA &RCGA &CLPSO &SaDE\\
    \midrule
    \multirow{5}{*}{Spherical}  &best &3.9321e-088	    &12.9783	&0.4381	&4.2719\\
                                  &median &2.7535e-085	&201.4712	&12.1488	&74.5293\\
                                  &worst &6.1007e-083	&701.7440	&402.9182	&606.4771\\
                                  &mean &3.3001e-084	&225.1812	&56.1705	&113.5754\\
                                  &st.dev. &1.1232e-083	&175.7136	&92.0918	&135.8662\\ \midrule

    \multirow{5}{*}{Rastrigin}  &best &0 &2.3910	&0.6742	&0.0943\\
                                  &median &0 &5.4255	&4.1572	&4.1060\\
                                  &worst &0 &11.4545	&9.5179	&11.4459\\
                                  &mean &0 &5.7737	&4.2650	&4.9949\\
                                  &st.dev. &0 &2.8549	&2.2645	&2.9296\\ \midrule

    \multirow{5}{*}{Griewank}   &best &0 &0.1613	&0.0607	&0.6381\\
                                  &median &0 &2.0181	&0.9858	&1.6781\\
                                  &worst &0 &7.3556	&2.3697	&6.6758\\
                                  &mean &0 &2.4742	&1.0333	&1.9131\\
                                  &st.dev. &0 &1.6757	&0.6281	&1.1905\\ \midrule

    \multirow{5}{*}{Rosenbrock} &best &0	&0.7506	&0.7522	&1.2169\\
                                  &median &0	&175.3620	&29.9583	&304.4219\\
                                  &worst &0	&8.5705e+03	&2.2800e+03	&7.4413e+03\\
                                  &mean &0	&1.0031e+03	&203.5033	&1.0847e+03\\
                                  &st.dev. &0	&2.1563e+03	&456.0105	&1.7735e+03\\ \midrule

    \multirow{5}{*}{Ackley}     &best &-8.8818e-016	&2.6351	&0.7704	&1.2860\\
                                  &median &-8.8818e-016	&8.0993	&4.5798	&8.6868\\
                                  &worst &-8.8818e-016	&15.1767	&10.8164	&13.6753\\
                                  &mean &-8.8818e-016	&8.6630	&5.0455	&8.4200\\
                                  &st.dev. &0	&3.2121	&2.6464	&3.2138\\ \midrule

    \multirow{5}{*}{Schaffer}   &best &0	    &0.0949	&0.0372	&0.0802\\
                                  &median &0	&0.2595	&0.1846	&0.2611\\
                                  &worst &0	&0.4233	&0.4182	&0.4451\\
                                  &mean &0	&0.2754	&0.2019	&0.2666\\
                                  &st.dev. &0	&0.1129	&0.1142	&0.0929\\ \midrule

    \multirow{5}{*}{Easom}      &best &-1 &-1.4112e-08	&-1.1388e-05	&-0.7633\\
                                  &median &-1 &-1.5910e-133	&-1.0708e-133	&-6.5793e-127\\
                                  &worst &-1 &0	&0	&0\\
                                  &mean &-1 &-9.9558e-10	&-3.7961e-07	&-0.0509\\
                                  &st.dev. &0	&3.5693e-09	&2.0792e-06	&0.1937\\ \midrule

    \multirow{5}{*}{Goldstein-Price}  &best &3 &3.6198	&3.0255	&3.1157\\
                                        &median &3 &12.7674	&15.0549	&16.5939\\
                                        &worst &3 &32.6912	&117.9224	&99.2768\\
                                        &mean &3 &15.1754	&20.4450	&22.0584\\
                                        &st.dev. &1.2148e-15 &9.5868	&22.9667	&20.8716\\
    \bottomrule
    \end{tabular}
\end{table*}

\begin{table*}[htbp]
\centering
\caption{Experimental results on test functions(10D)} \label{tab:results10D}
    \begin{tabular}{cccccc}
    \toprule
    Fcn                           &Statistic &MASTA &RCGA &CLPSO &SaDE\\
    \midrule
    \multirow{5}{*}{Spherical}  &best &1.6589e-026	&5.7655e+03	&499.6196	&3.1841e+03\\
                                  &median &7.0477e-025	&1.3771e+04	&1.8382e+03	&8.8483e+03\\
                                  &worst &4.4277e-023	&1.9348e+04	&9.6930e+03	&1.4828e+04\\
                                  &mean &3.7953e-024	&1.3597e+04	&2.4014e+03	&9.1901e+03\\
                                  &st.dev. &8.7862e-024	&3.5897e+03	&2.0599e+03	&2.8815e+03\\ \midrule

    \multirow{5}{*}{Rastrigin}  &best &0	&81.7097	&41.9288	&54.4110\\
                                  &median &0	&112.2439	&69.4195	&96.1864\\
                                  &worst &0	&138.2081	&107.2680	&125.4676\\
                                  &mean &0	&112.8742	&70.6708	&95.7393\\
                                  &st.dev. &0	&11.9848	&15.3291	&17.3010\\ \midrule

    \multirow{5}{*}{Griewank}   &best &0	&26.5595	&4.0602	&30.9557\\
                                  &median &0	&135.1928	&12.4689	&72.4932\\
                                  &worst &7.7716e-016	&161.5584	&62.9645	&155.9005\\
                                  &mean &4.4409e-017	&124.0792	&17.3946	&78.5525\\
                                  &st.dev. &1.4752e-016	&34.2075	&13.6905	&31.4477\\ \midrule

    \multirow{5}{*}{Rosenbrock} &best &2.5454e-010	&8.4705e+06	&2.6796e+04	&5.8034e+05\\
                                  &median &4.6860e-010	&2.2713e+07	&5.7970e+05	&8.1303e+06\\
                                  &worst &1.4606e-009	&4.5227e+07	&9.8850e+06	&6.7390e+07\\
                                  &mean &5.4748e-010	&2.5621e+07	&1.2024e+06	&1.3593e+07\\
                                  &st.dev. &2.5820e-010	&1.3481e+07	&1.9077e+06	&1.4995e+07\\ \midrule

    \multirow{5}{*}{Ackley}     &best &2.6645e-015	&17.4993	&7.6560	&15.9621\\
                                  &median &2.6645e-015	&19.6633	&12.9832	&18.6637\\
                                  &worst &6.2172e-015	&20.5070	&19.7384	&20.3793\\
                                  &mean &2.9014e-015	&19.6734	&13.1695	&18.4356\\
                                  &st.dev. &9.0135e-016	&0.5635	&3.1732	&1.2437\\
    \bottomrule
    \end{tabular}
\end{table*}

\begin{table*}[htbp]
\centering
\caption{Experimental results on test functions(20D)} \label{tab:results20D}
    \begin{tabular}{cccccc}
    \toprule
    Fcn                           &Statistic &MASTA &RCGA &CLPSO &SaDE\\
    \midrule
    \multirow{5}{*}{Spherical}  &best	&5.2958e-20	&2.4986e+04	&2.9887e+03	&3.2896e+03\\
                                    &median	&1.9835e-15	&3.6313e+04	&6.6088e+03	&2.4179e+04\\
                                    &worst	&2.4654e-14	&4.7295e+04	&1.1369e+04	&3.5268e+04\\
                                    &mean	&3.6743e-15	&3.6336e+04	&6.7626e+03	&2.3597e+04\\
                                    &st.dev.	&5.2892e-15	&7.5464e+03	&2.6801e+03	&7.5004e+03\\\midrule

    \multirow{5}{*}{Rastrigin}  &best	&0	&208.3856	&134.3816	&168.4955\\
                                    &median	&5.6843e-14	&265.1350	&179.8537	&245.3389\\
                                    &worst	&1.0658e-11	&297.2637	&223.1913	&304.6113\\
                                    &mean	&9.0286e-13	&262.7806	&175.5388	&243.7538\\
                                    &st. dev.	&2.1777e-12	&23.4442	&22.9825 &28.8339 \\ \midrule

    \multirow{5}{*}{Griewank}   &best	&0	&220.6985	&20.7024	&41.2874\\
                                    &median	&2.6090e-15	&329.1324	&51.6339	&191.1037\\
                                    &worst	&2.6423e-14	&446.4131	&113.4861	&355.0719\\
                                    &mean	&6.1395e-15	&339.4656	&56.3813	&186.5737\\
                                    &st. dev.	&7.6550e-15	&61.0506	&23.5577 &57.0128\\ \midrule

    \multirow{5}{*}{Rosenbrock} &best	&8.9275e-11	&5.7734e+07	&3.8928e+05	&1.9639e+07\\
                                    &median	&5.1014e-10	&1.1769e+08	&3.2493e+06	&4.1686e+07\\
                                    &worst	&4.0037e-06	&1.7354e+08	&8.1657e+07	&1.1972e+08\\
                                    &mean	&1.3393e-07	&1.7354e+08	&7.0817e+06	&4.8172e+07\\
                                    &st. dev.	&7.3089e-07	&2.8629e+07	&1.4577e+07	&2.5965e+07 \\ \midrule

    \multirow{5}{*}{Ackley}     &best	&6.2172e-15	&19.9040	&10.8397	&17.3019\\
                                    &median	&1.3323e-14	&20.3015	&15.9606	&19.5759\\
                                    &worst	&2.0099e-12	&20.7213	&17.9238	&20.4883\\
                                    &mean	&1.4264e-13	&20.3455	&15.2953	&19.3847\\
                                    &st. dev.	&4.2496e-13	&0.2478	&1.7842	&0.9224\\
    \bottomrule
    \end{tabular}
\end{table*}

\begin{table*}[htbp]
\centering
\caption{Experimental results on test functions(30D)} \label{tab:results30D}
    \begin{tabular}{cccccc}
    \toprule
    Fcn                           &Statistic &MASTA &RCGA &CLPSO &SaDE\\
    \midrule
    \multirow{5}{*}{Spherical}  &best	&9.9298e-17	&5.0224e+04	&5.9480e+03	&1.7296e+04\\
                                    &median	&2.2586e-15	&5.8200e+04	&1.1819e+04	&4.1497e+04\\
                                    &worst	&2.9259e-13	&7.4783e+04	&1.8616e+04	&6.5872e+04\\
                                    &mean	&5.4250e-14	&5.9707e+04	&1.1357e+04	&4.2381e+04\\
                                    &st. dev.	&8.6142e-14	&7.3451e+03	&3.5588e+03	&1.3070e+04\\\midrule

    \multirow{5}{*}{Rastrigin}  &best	    &0	&357.3603	&238.6496	&290.5916\\
                                    &median	&4.5475e-13	&427.7500	&296.7514	&368.8724\\
                                    &worst	&3.1434e-10	&489.4174	&418.8757	&434.9728\\
                                    &mean	&1.9376e-11	&421.9971	&299.6580	&369.5383\\
                                    &st. dev.	&5.8378e-11	&29.3700	&37.7247	&33.3396\\ \midrule

    \multirow{5}{*}{Griewank}   &best	&0	&439.6893	&57.1764	&201.6535\\
                                    &median	&1.4433e-15	&563.4664	&57.1764	&384.5188\\
                                    &worst	&4.8583e-13	&700.3678	&258.4568	&590.8669\\
                                    &mean	&5.8387e-14	&562.2621	&113.3711	&384.1066\\
                                    &st. dev.	&1.1433e-13    	&71.9820	&45.9129	&94.2135\\ \midrule

    \multirow{5}{*}{Rosenbrock} &best	&1.5584e-10	&1.5442e+08	&1.7075e+06	&1.9372e+07\\
                                    &median	&5.9105e-10	&2.3835e+08	&9.0885e+06	&9.7060e+07\\
                                    &worst	&1.0345e-08	&3.4958e+08	&1.3413e+08	&2.5800e+08\\
                                    &mean	&9.3925e-10	&2.3523e+08	&1.4109e+07	&9.8533e+07\\
                                    &st. dev.	&1.7887e-09	&4.4129e+07	&2.3406e+07	&5.7739e+07\\ \midrule

    \multirow{5}{*}{Ackley}     &best	&2.0517e-13	&20.2510	&13.8995	&16.8680\\
                                    &median	&7.6117e-13	&20.7147	&16.0119	&19.6832\\
                                    &worst	&7.2502e-12	&20.9452	&17.7340	&20.5892\\
                                    &mean	&9.6877e-13	&20.6411	&15.9828	&19.4469\\
                                    &st. dev.	&1.2632e-12	&0.1995	&1.0368	&0.9342\\
    \bottomrule
    \end{tabular}
\end{table*}

\begin{table*}[htbp]
\centering
\caption{Experimental results on test functions(50D)} \label{tab:results50D}
    \begin{tabular}{cccccc}
    \toprule
    Fcn                           &Statistic &MASTA &RCGA &CLPSO &SaDE\\
    \midrule
    \multirow{5}{*}{Spherical}  &best	&8.5786e-14	&1.0249e+05	&1.7001e+04	&3.7394e+04\\
                                    &median	&2.2334e-13	&1.1797e+05	&2.5192e+04	&6.8738e+04\\
                                    &worst	&4.5530e-12	&1.3835e+05	&3.6831e+04	&9.7637e+04\\
                                    &mean	&6.3530e-13	&1.1973e+05	&2.6623e+04	&6.7411e+04\\
                                    &st. dev.	&9.6000e-13	&8.8766e+03	&5.1924e+03	&1.5539e+04\\\midrule

    \multirow{5}{*}{Rastrigin}  &best	&2.2737e-13	&676.1611	&450.7825	&519.3628\\
                                    &median	&3.5811e-12	&759.2857	&506.6920	&660.9823\\
                                    &worst	&1.6800e-09	&834.4782	&572.7979	&723.2395\\
                                    &mean	&7.0215e-11	&758.3887	&504.8464	&640.5101\\
                                    &st. dev.	&3.0526e-10	&33.9045	&31.8814	&58.0821\\ \midrule

    \multirow{5}{*}{Griewank}   &best	&1.3212e-14	&845.0041	&152.0123	&343.0254\\
                                    &median	&3.3662e-13	&1.0764e+03	&239.6336	&597.2199\\
                                    &worst	&6.0544e-12	&1.2260e+03	&407.3338	&914.6162\\
                                    &mean	&9.5671e-13	&1.0544e+03	&237.4946	&605.9091\\
                                    &st. dev.	&1.4539e-12	&109.5325	&58.0561	&145.1355\\ \midrule

    \multirow{5}{*}{Rosenbrock} &best	&4.6608e-10	&3.1104e+08	&1.2155e+07	&5.2011e+07\\
                                    &median	&9.2241e-10	&4.5546e+08	&2.8330e+07	&1.7362e+08\\
                                    &worst	&8.0216e-08	&5.7721e+08	&6.4922e+07	&5.2442e+08\\
                                    &mean	&6.3139e-09	&4.4997e+08	&3.0119e+07	&2.0412e+08\\
                                    &st. dev.	&1.6969e-08	&6.7224e+07	&1.2155e+07	&1.0509e+08\\ \midrule

    \multirow{5}{*}{Ackley}     &best	&2.1734e-12	&20.5232	&15.3432	&18.0430\\
                                    &median	&1.0515e-11	&20.8280	&17.0050	&20.0115\\
                                    &worst	&3.3675e-11	&21.0206	&18.3675	&20.7485\\
                                    &mean	&1.2473e-11	&20.8079	&17.0503	&19.8433\\
                                    &st. dev.	&7.7607e-12	&0.1228	&0.7975	&0.5768\\
    \bottomrule
    \end{tabular}
\end{table*}

\subsection{Test functions}
In order to test the performance of MASTA, eight common benchmark functions (i.e., the Spherical function, Rastrigin function, Griewank function, Rosenbrock function, Ackley function, Schaffer function, Easom function and Goldstein-Price function) are selected for the experiment. Among them, these functions are multidimensional functions of various modals and the other three are two dimensional functions, which are listed in Table \ref{tab:BenchmarkFunctions}.
Moreover, the performance of MASTA is further compared with three state-of-the-art optimization algorithms, that is, RCGA \cite{Tran2010}, CLPSO\cite{liang2006comprehensive}, and SaDE \cite{Wang2011} in the experiments.

\begin{figure*}[htbp]
\centering

\subfigure[Spherical]{
\begin{minipage}{6cm}
\centering
\includegraphics[width=6cm]{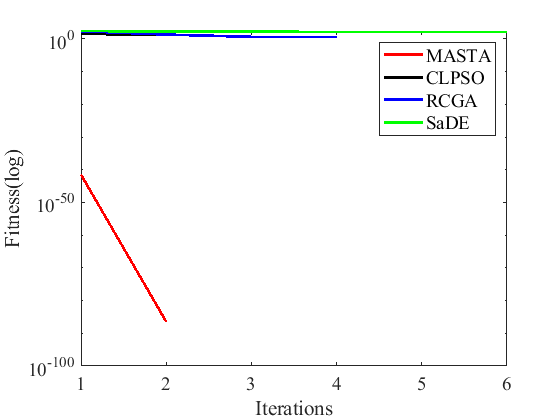}
%\caption{spherical}
\end{minipage}%
}%
\subfigure[Rosenbrock]{
\begin{minipage}{6cm}
\centering
\includegraphics[width=6cm]{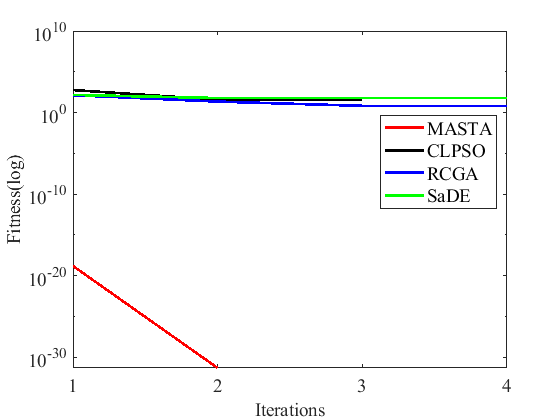}
%\caption{rosenbrock}
\end{minipage}
}

\subfigure[Rastrigin]{
\begin{minipage}{5.5cm}
\centering
\includegraphics[width=5.5cm]{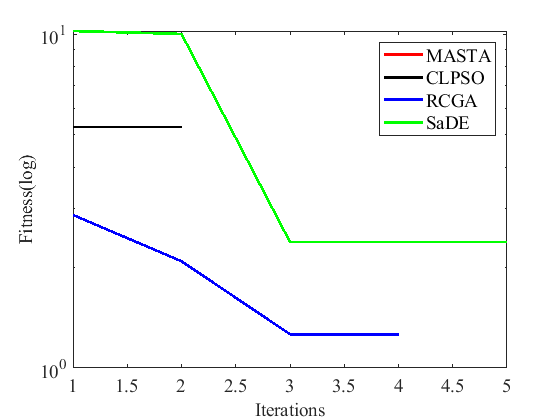}
%\caption{rastrigin}
\end{minipage}
}
\subfigure[Griewank]{
\begin{minipage}{5.5cm}
\centering
\includegraphics[width=5.5cm]{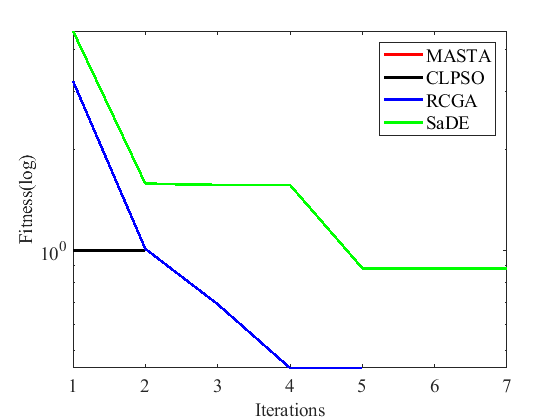}
%\caption{griewank}
\end{minipage}%
}%

\subfigure[Ackley]{
\begin{minipage}{6cm}
\centering
\includegraphics[width=6cm]{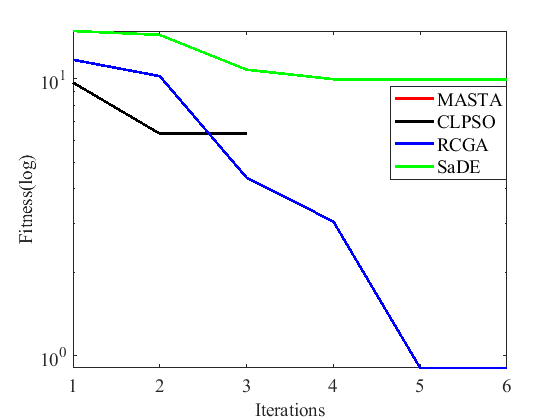}
%\caption{Ackley}
\end{minipage}%
}%
\subfigure[Schaffer]{
\begin{minipage}{6cm}
\centering
\includegraphics[width=6cm]{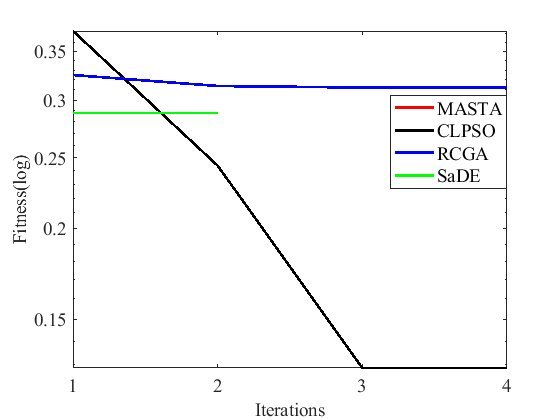}
%\caption{schaffer}
\end{minipage}
}

\subfigure[Easom]{
\begin{minipage}{6cm}
\centering
\includegraphics[width=6cm]{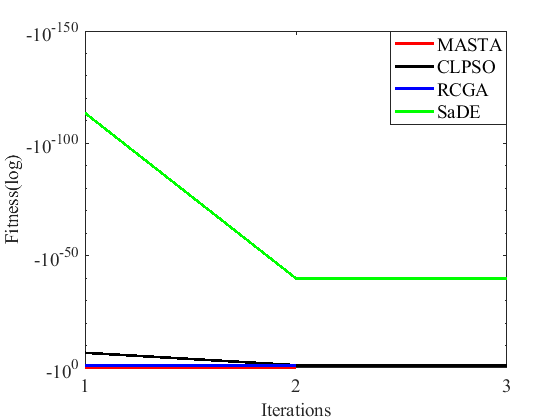}
%\caption{easom}
\end{minipage}%
}%
\subfigure[Goldstein-Price]{
\begin{minipage}{6cm}
\centering
\includegraphics[width=6cm]{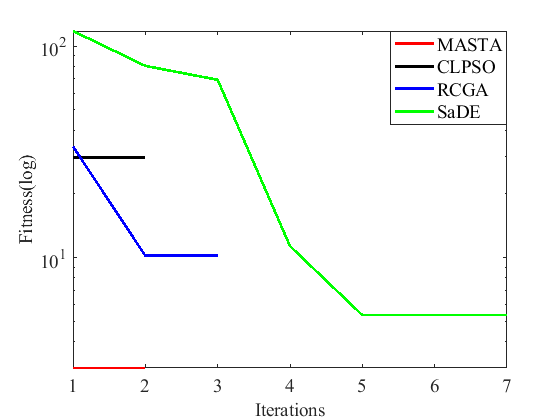}
%\caption{goldstein-price}
\end{minipage}
}

\caption{Average fitness of the two-dimensional functions (2D)}
\label{fig:average_fitness_2D}
\end{figure*}

\begin{figure*}[htbp]
\centering

\subfigure[Spherical]{
\begin{minipage}{6cm}
\centering
\includegraphics[width=6cm]{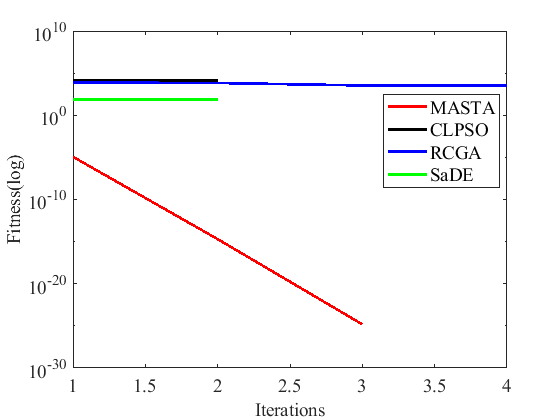}
%\caption{Spherical}
\end{minipage}%
}%
\subfigure[Rosenbrock]{
\begin{minipage}{6cm}
\centering
\includegraphics[width=6cm]{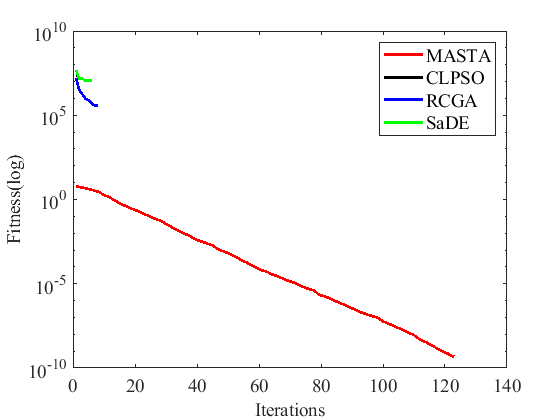}
%\caption{Rosenbrock}
\end{minipage}
}

\subfigure[Rastrigin]{
\begin{minipage}{5.6cm}
\centering
\includegraphics[width=5.6cm]{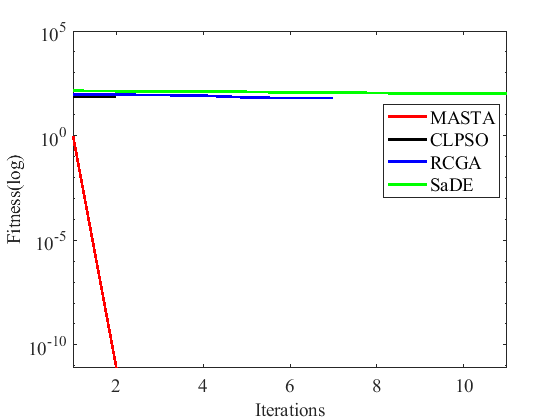}
%\caption{Rastrigin}
\end{minipage}
}
\subfigure[Griewank]{
\begin{minipage}{5.6cm}
\centering
\includegraphics[width=5.6cm]{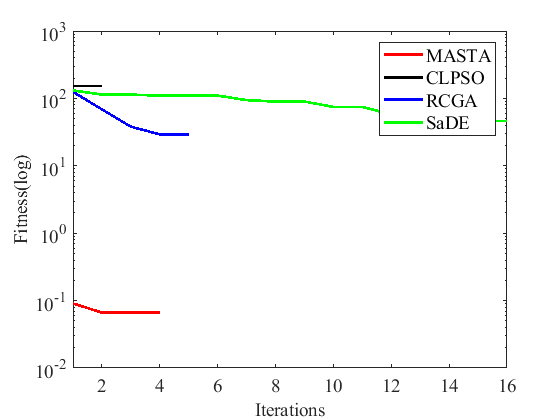}
%\caption{Griewank}
\end{minipage}%
}%

\subfigure[Ackley]{
\begin{minipage}{6cm}
\centering
\includegraphics[width=6cm]{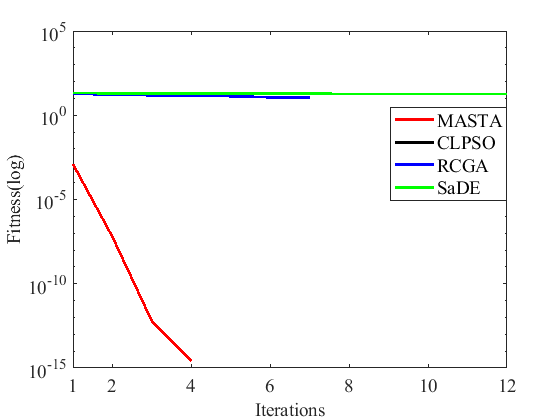}
%\caption{Ackley}
\end{minipage}%
}%

\caption{Average fitness of the ten-dimensional functions (10D)}
\label{fig:average_fitness_10D}
\end{figure*}

\begin{figure*}[htbp]
\centering

\subfigure[Spherical]{
\begin{minipage}{6cm}
\centering
\includegraphics[width=6cm]{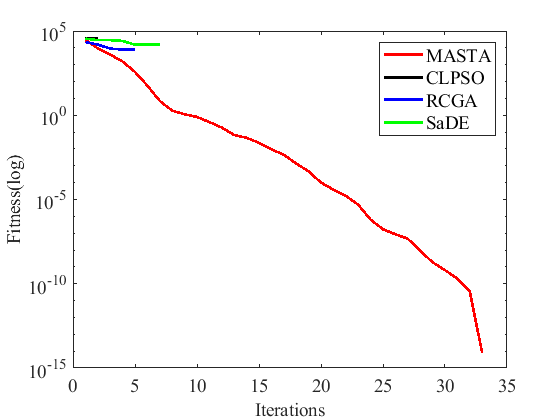}
%\caption{Spherical}
\end{minipage}%
}%
\subfigure[Rosenbrock]{
\begin{minipage}{6cm}
\centering
\includegraphics[width=6cm]{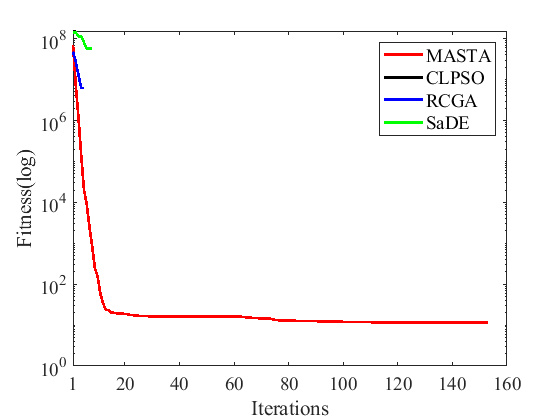}
%\caption{Rosenbrock}
\end{minipage}
}

\subfigure[Rastrigin]{
\begin{minipage}{5.6cm}
\centering
\includegraphics[width=5.6cm]{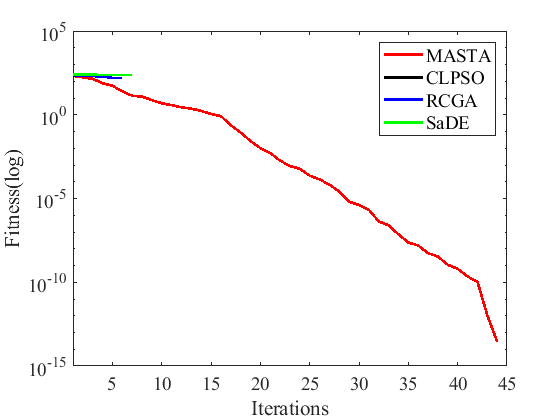}
%\caption{Rastrigin}
\end{minipage}
}
\subfigure[Griewank]{
\begin{minipage}{5.6cm}
\centering
\includegraphics[width=5.6cm]{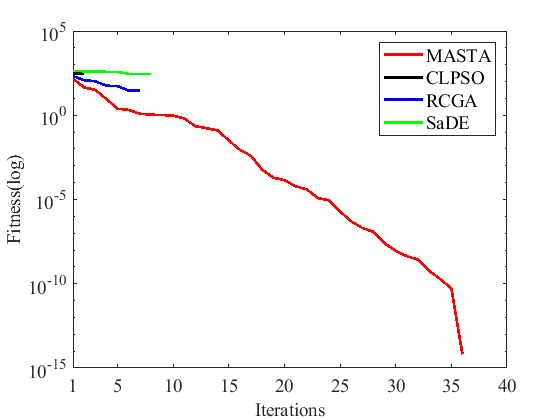}
%\caption{Griewank}
\end{minipage}%
}%

\subfigure[Ackley]{
\begin{minipage}{6cm}
\centering
\includegraphics[width=6cm]{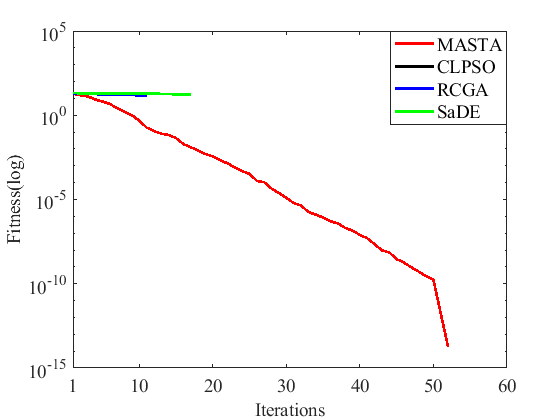}
%\caption{Ackley}
\end{minipage}%
}

\caption{Average fitness of the ten-dimensional functions (20D)}
\label{fig:average_fitness_20D}
\end{figure*}

\begin{figure*}[htbp]
\centering

\subfigure[Spherical]{
\begin{minipage}{6cm}
\centering
\includegraphics[width=6cm]{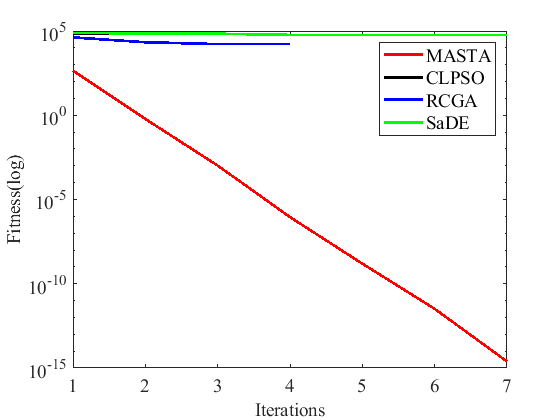}
%\caption{Spherical}
\end{minipage}%
}%
\subfigure[Rosenbrock]{
\begin{minipage}{6cm}
\centering
\includegraphics[width=6cm]{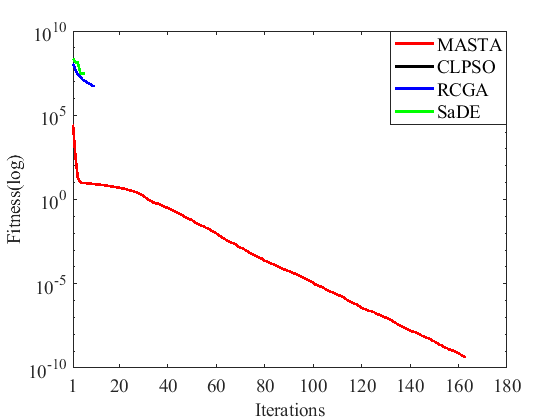}
%\caption{Rosenbrock}
\end{minipage}
}

\subfigure[Rastrigin]{
\begin{minipage}{5.6cm}
\centering
\includegraphics[width=5.6cm]{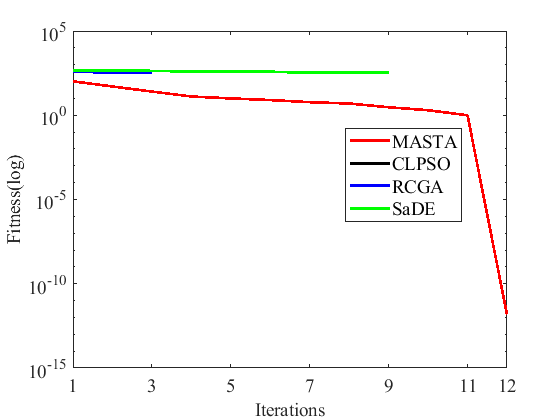}
%\caption{Rastrigin}
\end{minipage}
}
\subfigure[Griewank]{
\begin{minipage}{5.6cm}
\centering
\includegraphics[width=5.6cm]{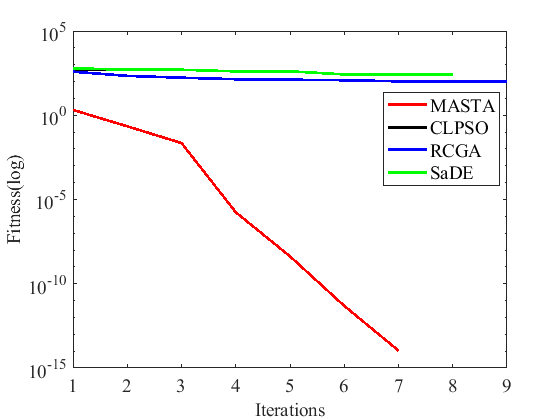}
%\caption{Griewank}
\end{minipage}%
}%

\subfigure[Ackley]{
\begin{minipage}{6cm}
\centering
\includegraphics[width=6cm]{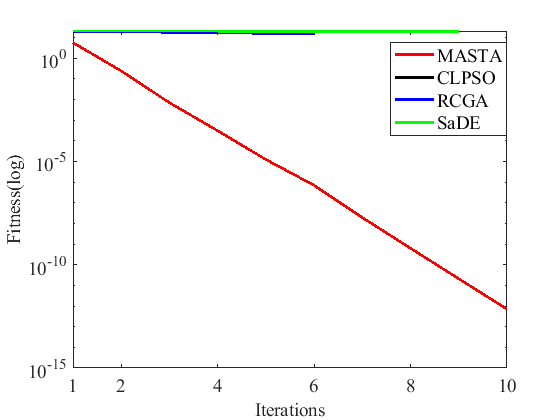}
%\caption{Ackley}
\end{minipage}%
}

\caption{Average fitness of the ten-dimensional functions (30D)}
\label{fig:average_fitness_30D}
\end{figure*}

\begin{figure*}[htbp]
\centering

\subfigure[Spherical]{
\begin{minipage}{6cm}
\centering
\includegraphics[width=6cm]{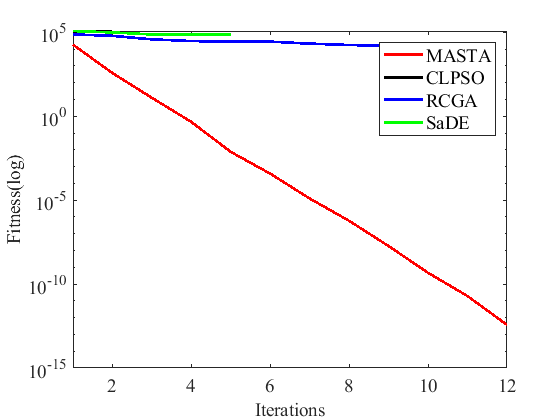}
%\caption{Spherical}
\end{minipage}%
}%
\subfigure[Rosenbrock]{
\begin{minipage}{6cm}
\centering
\includegraphics[width=6cm]{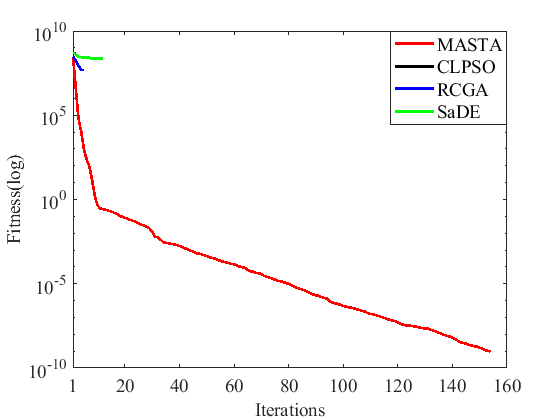}
%\caption{Rosenbrock}
\end{minipage}
}

\subfigure[Rastrigin]{
\begin{minipage}{5.6cm}
\centering
\includegraphics[width=5.6cm]{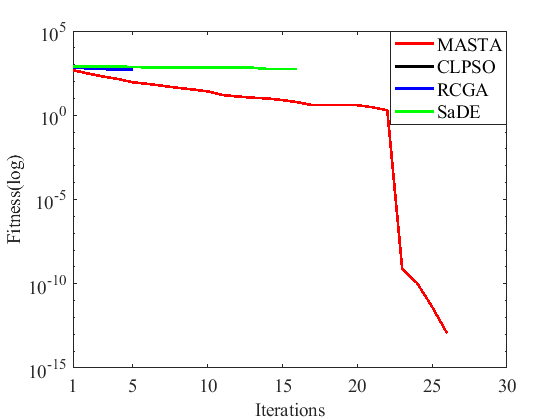}
%\caption{Rastrigin}
\end{minipage}
}
\subfigure[Griewank]{
\begin{minipage}{5.6cm}
\centering
\includegraphics[width=5.6cm]{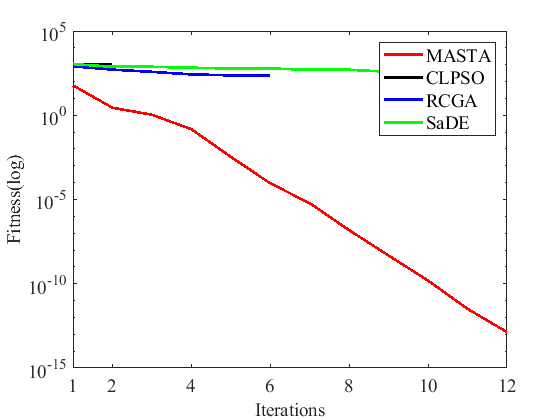}
%\caption{Griewank}
\end{minipage}%
}%

\subfigure[Ackley]{
\begin{minipage}{6cm}
\centering
\includegraphics[width=6cm]{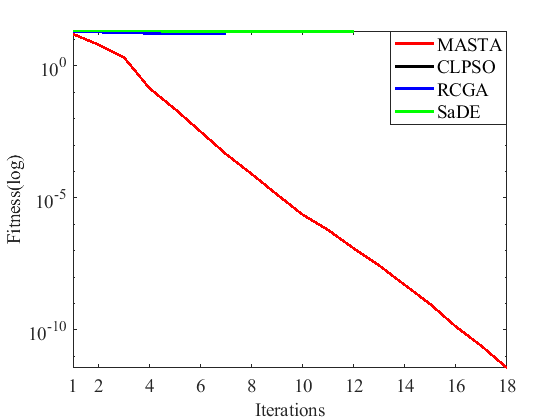}
%\caption{Ackley}
\end{minipage}%
}

\caption{Average fitness of the ten-dimensional functions (50D)}
\label{fig:average_fitness_50D}
\end{figure*}

\subsection{Parameters setting}
All the experiments of the proposed MASTA algorithm are performed in MATLAB R2019a on Windows 10 64 bit with a Intel Core i5 2.5GHz CPU and 6.0GB RAM.
The parameters in the proposed MASTA are set as follows: the number of states is $SE=20$, the communication frequency is $CF=50$, all the control parameters of transformation operators are set as 1. In this paper, the exponential way is accepted for its rapidity, of which the base is 2 in the experiment.
As for RCGA, we use the same parameter settings as in \cite{Tran2010}. Then, for CLPSO and SaDE, we use the MATLAB codes provided by the author in \cite{liang2006comprehensive,Wang2011} with minor revisions for this experiment.
Commonly, the variation of a parameter follows a linear, exponential or a logistic way.
Programs were performed for 30 runs independently and the population scale is 30 for each run in experiments.

\subsection{Results and discussion}
For comparison, some common statistics are introduced. The best means the minimum of the results, the worst indicates the maximum of the results, and then it follows the mean, median and standard deviation (st.dev.). In some way, these statistics are able to evaluate the search ability and solution accuracy, reliability and convergence as well as stability. To be more specific, the best indicates the global search ability and solution accuracy, the worst and the mean signify the reliability and convergence, while the median and standard deviation correspond to the stability.

Results for two dimensional functions optimization are listed in Table \ref{tab:results2D} and results for ten dimensional functions optimization can be found in Table \ref{tab:results10D}.
On the other hand, illustrations of the average fitness in 30 simulations are given in Figure \ref{fig:average_fitness_2D}, and Figure \ref{fig:average_fitness_10D}, respectively.
The average fitness curve can visually depict the search ability and convergence performance. In the following the analysis of the results for each function is discussed separately.

For the Spherical function, from the experimental results, we can see that all of the algorithms can
find the global optimum with high solution precision and have good reliability as well as stability for this function in terms of two and ten dimensions.
For the Rastrigin function, as can be seen from the experimental results, all of the algorithms are able to find the global optimum with high solution precision and have good reliability as well as stability for all the two dimensional function. Moreover, for the ten dimensional function, the global optimum can also be found by all algorithms; however, the MASTA algorithm has much better statistical performances than other three algorithms especially described by the worst. RCGA and SaDE can not achieve the best occasionally, and the mean of RCGA and SaDE are not satisfactory.
For the Griewank function, from the experimental results, we can find that, for the two dimensional function, all the three algorithms, i.e., MASTA, CLPSO and SaDE, have the very good ability to achieve the best solutions and have both very good reliability and stability while the RCGA has the worst reliability and stability. For the ten dimensional function, MASTA has the very excellent ability to achieve the best solutions with both very good reliability and stability, and the other three algorithms show very poor reliability and stability.
For the Rosenbrock function, we can see from the experimental results, for the dimensional function, all the algorithms can find the global optimum with very good reliability and stability for this function in space, while the RCGA shows very poor reliability and stability for the function. For the ten dimensional problem, only MASTA is able to find the optimal solutions with very good reliability and stability while the other three algorithms RCGA, CLPSO and SaDE show very poor reliability and stability for the function.
For the Ackley function, as can be seen from the experimental results, all these algorithms can find the global optimum with high solution precision and have very good reliability as well as stability for all the functions in terms of two and ten dimensions, and because the st.dev. approaches zero the statistical performances of results are very satisfactory for all these algorithms.
For the Schaffer function, for the two dimensional function, we can see that the MASTA algorithm can find the global optimum with very good reliability as well as stability for the function, while the other three algorithms RCGA, CLPSO and SaDE show very poor reliability and stability for the function.
For both the Easom and Goldstein-Price functions, we can see that all these algorithms can find the global optimum with very high reliability as well as stability for these functions.
In addition, from the Figure \ref{fig:average_fitness_2D}, and Figure \ref{fig:average_fitness_10D}, we can also see that MASTA algorithm can converge to the global optimum in a very fast rate for all these functions in the experiments.
Thus, according to the experimental results, we can see that the MASTA algorithm can find the global optimum with high solution precision and have good reliability as well as stability for all these two dimensional and ten dimensional functions. Moreover, the searching time required for MASTA is infinity in theory, which is the consequence of random search methodology. However, in practice, we can stop the iteration process by presetting some criteria, for example, the prescribed maximum iterations, or when the fitness is unchanged for a number of times, and we use the no improvement of the found global optimal solution in the MASTA algorithm as the iteration criterion in the paper.

\section{Conclusion}
In this paper, we propose an effective multiagent-based state transition algorithm named MASTA. In the proposed MASTA algorithm, we apply the state transition algorithm to the population of solutions and perform an effective communication operation to update the population in a iterative process to find the best optimal solution from the search space effectively and efficiently. Based on the experiments on some common benchmark functions with different dimensions, the results have demonstrated that the MASTA algorithm can successfully find the best optimal solutions in all these complex functions. In addition, compared with some state-of-the-art optimization algorithms, the MASTA algorithm has shown very superior and comparable global search performance.

%\section*{Acknowledgment}
%% Can use something like this to put references on a page
%% by themselves when using endfloat and the captionsoff option.
%\ifCLASSOPTIONcaptionsoff
%  \newpage
%\fi



%\bibliographystyle{IEEEtran}
%\begin{thebibliography}{1}
%\end{thebibliography}

\bibliographystyle{ieeetr}
\bibliography{tzy}

%\begin{biography}[{\includegraphics[width=1in,height=1.25in,clip,keepaspectratio]{mshell}}]{Michael Shell}
%\end{IEEEbiography}
%% insert where needed to balance the two columns on the last page with
%%\newpage


%\vfill


\end{document}